\documentclass[11pt]{llncs}

 \usepackage{graphics} 
\usepackage{epsfig} 
\usepackage{caption}
\usepackage{float}

\begin{document}


\title{A Concise Proof of Discrete Jordan Curve Theorem}
\author{Li Chen}

\institute{University of the District of Columbia}

\maketitle

\begin{abstract}

In geometry and graph theory, the Jordan curve theorem has both educational and academic roles. There have been many
proofs published in educational  journals in mathematics， notably ``American Mathematical Monthly.''\cite{Mae,Tho,FW,Hal}
This paper gives a concise proof of the Jordan curve theorem on discrete surfaces.
We also embed the discrete surface in the 2D plane to prove the original version
of the Jordan curve theorem. This paper is a simple version of \cite{Che13}.
We seek to clarify and simplify some statements and proofs. Again, the purpose
of this paper is to make the proof of the theorems easier to understand.
We omit some definitions of concepts, which can be found in \cite{Che04,Che12,Che99,Che13}. In revision 2, we added another appendix
to make a self-contained proof on verifying simple connectedness of the Euclidean plane in this paper.
In this revision, we added a special case for the proof of Theorem 3 in Appendix B that was found when we were revising a new paper for high dimensional contraction
\cite{Che17}. It was easy to resolve in 2D. We put it in Appendix C of this paper.
\footnote {Readers who do not want to know the details of the concepts in discrete surfaces can begin
in Section 2 to get an idea of the proofs based on the intuitive meaning of the concepts used
in this paper.}
\end{abstract}



\section{What is the Jordan Curve Theorem?}

In this paper, we give a straightforward proof of the Jordan curve theorem in 2D discrete spaces
with respect to the general definition of discrete curves, surfaces, and manifolds discussed in ~\cite{Che97,Che04}.
The Jordan curve theorem states that
a simple and closed  curve separates a simply connected surface into
two components.  Based on the definition of discrete surfaces,
we give three reasonable definitions of simply connected spaces in discrete spaces.
Theoretically, these three definitions are equivalent.

For the Jordan curve theorem, O. Veblen in 1905 wrote a  paper~\cite{Veb} that
was regraded as the first correct proof of this fundamental theorem on 2D Euclidean plane.
The first discrete proof was given by W.T. Tutte on planar graphs in 1979 ~\cite{Tut}.
Recently, researchers still show
considerable interests in the Jordan curve theorem using formalized proofs in  computers ~\cite{Duf}.

In 1999, L. Chen attempted to prove  the discrete Jordan curve theorem for 2D discrete manifolds
without using 2D Euclidean space ~\cite{Che99}. Chen added some missing parts of this proofs in 2013 \cite{Che13}.
In \cite{Che13}, Chen adopted some original ideas from  Veblen's paper and gave a proof of
this theorem in discrete form.

In this paper, we want to give a simple and easy version of the proofs given in \cite{Che13} without
gone through many definitions of discrete geometry.

A discrete surface can be viewed as a digital surface or a triangulated 2D meshes. Intuitively,
a discrete 2-cell in this paper is a smallest or minimal unite for 2D objects. A discrete 2-cell
cannot be slitted into two other 2-cells. The union of two (discrete) 2-cells ,$A$ and $B$ will not be a two cell if
$A\neq B$. (In discrete space, an object is mainly formed by vertices or points. The others such as edges, faces
are based on human's interpretation.)
In other words, the smallest unite does not contain any other 2-cell. When we say a 2-cell, we need to think this ``thing'' as
the smallest one containing some vertices and edges. Please note that this definition is different from
the standard definition in topology.

We know that a discrete surface can be naturally embedded to Euclidean plane or a closed continuous surface such as
a sphere. Indeed, a discrete surface can be
easily embedded to a 3D or higher dimensional Euclidean spaces.

Let us first review some concepts of discrete curves in \cite{Che99,Che04}: (1) A simple path
is called a (discrete) pseudo-curve, (2) a simple semi-curve
can be a curve or a surface-cell (2-cell) (in discrete form, a set containing all vertices
means that the set contains all 2-cell in this set. At least it is true for this paper. See
details in \cite{Che04} for general non simply connected surfaces.)
(3) A simple curve must not contain any
proper subset that is a 2-cell.

It is obvious, if we define a simple path (a pseudo-curve in this paper) a discrete curve, there is
no Jordan Curve theorem in discrete space. This is because that the inner part of a 2-cell
is empty in graph-structures. This little difference will not affect the continuous version
of the theorem. In discrete cases, we will make a point called ``the Veblen point''
to deal with this when we desire a closer version of the theorem to the continuous
case.

We defined regular points in a discrete manifolds in \cite{Che04}. It just means that
every point on surface has a neighborhood that is similar to a 2D disk (homeomorphic equivalence in
topology).

Let $S(p)$ be the neighborhood of $p$ in $S$. If $S$ is a triangulated surface, then

\[ S(p) =\{v | v ,  p \mbox{  are adjacent.} \}\cup \{p\}\]

\noindent or $S(p) =\{v | d(v,p)\le 1\}$ where $d$ is graph-distance.  For any type of discrete
surfaces, $ S(p) =\{v | v  \mbox{ is in a 2-cell that contains} p \}\cup \{p\}$.

We proved that for a  discrete surface $S$, if $p$ is a inner and regular
point of $S$, then there exists a simple cycle containing
all points in $S(p)-\{p\}$ in $S$.
This result is particularly important in our proof, but it is very intuitive too.
We present this result as Lemma 1 \cite{Che04,Che13}.

\begin{lemma}
For a  discrete surface $S$, if $p$ is a inner and regular
point of $S$, then there exists a simple cycle containing
all points in $S(p)-\{p\}$ in $S$ where $S(p)$ is the neighborhood of $p$ in $S$.
\end{lemma}

In topology, the formal description of the Jordan curve theorem is:
A simply closed curve $J$ in a plane $\Pi$ decomposes $\Pi - J$ into two
components.~\cite{Lef,New}  In fact, this theorem holds for any simply connected 2D surface.
A plane is a simply connected surface in Euclidean space, but this theorem is not true
for a general continuous surface. For example, the boundary of a donut.

We now introduce the meaning of discrete deformation and simply connected discrete surfaces.

What is a simply connected continuous surface?
A connected topological space $T$ is simply connected
if for any point $p$ in $T$, any simply closed curve containing $p$ can be
contracted to $p$. The contraction is a continuous mapping among
a series of closed continuous curves. ~\cite{New}
So, we first need the concept of ``discrete contraction.''

Here we also try to make the definition of simply connected discrete surfaces to be simpler.
The more detailed definition was in \cite{Che12,Che04}.

The ``contraction'' means the sequence of curves will be getting smaller and smaller until
they shrink to a point. These curves do not cross each other. (They do not have to be cross each other.
This is a key too.)

So our concept of contraction is that a discrete curve will be graduate varied to a ``smaller'' one,
and keep the process until ``shrink'' to a point. These discrete curves do not cross each other.
they can share some points.

In order to keep the concepts simple to understand, we defined the
gradual variation between two simple paths in \cite{Che04,Che12,Che13}.
We defined discrete deformation
among discrete pseudo curves. And finally, we define the contraction of curves
is a type of discrete deformation. See \cite{Che04,Che12} for more
details of the definitions. The reader can just use the natural interpretation of
definitions.  In paper, we assume the discrete surface is both regular and orientable too.
The algorithm to decide if a discrete surface is orientable can be found in  \cite{Che04}.

Intuitively, two simple paths $C$ and $C'$ are said to be gradually varied if there is no hole
in between $C$ and $C'$. In addition, there is no jump from $C$ to $C'$.

More formally, two simple paths $C$ and $C'$ are called {\bf gradually varied} if $C \cup C'$ consists of
1-cells and 2-cells where no cycle in $C \cup C'$ that is not a 2-cell or the union of 2-cells.
In other words, Assume $E(C)$ denotes all edges in
path $C$.  Let $XorSum(C,C') = (E(C)-E(C'))\cup (E(C')-E(C))$.
$XorSum$  is called $sum ( modulo 2)$ in  Newman's book \cite{New}.
$C$ and $C'$ are gradually varied iff $XorSum(C,C')$ are the union of 2-cells.

To prove the Jordan curve theorem, we need to describe what the
disconnected components are by means of separated
from a simple curve $C$? It means that any path from a component
to another must include at least a point in $C$.
It also means that this linking path must cross-over the curve $C$.
We will define the concept of ``cross-over'' in the following.

Because a surface-cell $A$ is a closed path, we can define two
orientations (normals ) to $A$: clockwise and counter-clockwise.
Usually, the orientation of a 2-cell is not a critical issue.
However, for the proof of the Jordan curve theorem it is necessary.

In other words, a pseudo-curve which is a set of points has no ``direction,''
but as a path $P=\{p_0,\cdots,p_n\}$, it has its own ``travel direction'' from $p_0$ to $p_n$.

For two paths $C$ and $C'$,
which are gradually varied, if a 2-cell $A$
is in $G(C\cup C')$, the orientation of $A$ with respect to $C$ is determined by
the first pair of points $(p,q) \in C \cap A$ and $C= ... p q ... $ .
Moreover, if a 1-cell of $A$ is in $C$, then the orientation of $A$ is
fixed with respect to $C$.

According to Lemma 1, $S(p)$ contains all adjacent points of $p$ and
$S(p)-\{p\}$ is a simple cycle---there is a cycle containing all points in $S(p)-\{p\}$.

We assume that cycle $S(p)-\{p\}$ is always oriented
clockwise. For two points $a, b \in S(p)-\{p\}$,
there are two simple cycles containing the path $a\rightarrow p \rightarrow b$ :
(1) a cycle from $a$ to $p$ to $b$ then moving clockwise to $a$, and
(2) a cycle from $a$ to $p$ to $b$ then moving counter-clockwise to $a$.
See Fig. 1(a).

 It is easy to see that the simple cycle $S(p)-\{p\}$ separates
$S-\{S(p)-\{p\}\}$ into at least two connected components because from
$p$ to any other points in $S$ the path must contain a point in $S(p)-\{p\}$.
$S(p)-\{p\}$ is an example the Jordan curve.

\begin{figure}[h]
	\begin{center}

   \epsfxsize=4.5in
   \epsfbox{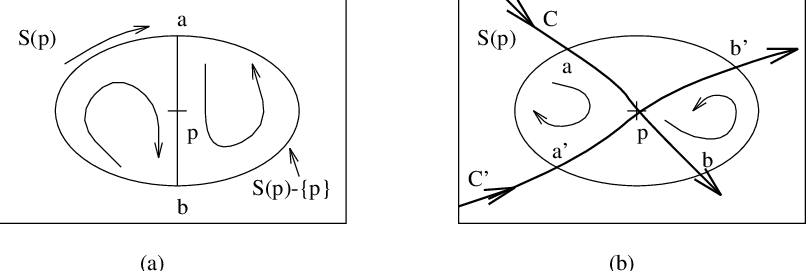}   

\caption{$S(p)$ and Cross-over at $p$: (a) Two adjacent points $a$ and $b$ of $p$ in
 $S(p)$, and (b) an example for two cross-over paths.}

\end{center}
\end{figure}

\begin{definition}
Two simple paths $C$ and $C'$ are said to be  ``cross-over'' each other if
there are points $p$ and $q$ ($p$ may be the same as $q$) such that
$C = ... a p b...s q t...$ and $C'= ...a' p b...s q t'...$ where
$a\ne a'$ and $t\ne t'$. The cycle
$a p a' ... a$ without $b$ in $S(p)$ and the cycle $q t ...t' q$ without
$s$ in $S(q)$ have different orientations with respect to $C$.
\end{definition}

For example, in Fig. 1 (b), $C$ and $C'$ are ``cross-over'' each other.
When  $C$ and $C'$ are not ``cross-over'' each other, we will say
that $C$ is at a side of $C'$.  We also say that $C$ and $C'$ in the above Lemma  are side-gradually
varied.

\begin{lemma}
If two simple paths $C$ and $C'$ are not cross-over each other, and they
are gradually varied, then every surface-cell in $G(C\cup C')$ has the
same orientation with respect to the ``travel direction'' of $C$ and
opposite to the ``travel direction'' of $C'$.
\end{lemma}

Intuitively, a simply connected set is such a set so that for any point,
every simple cycle containing
this point can contract to the point.

\begin{definition}
A simple cycle $C$ can contract to a point $p\in C$
if there exist a series of simple cycle, $C=C_0, ..., p=C_{n}$:
(1) $C_i$ contains $p$ for all $i$;
(2) If $q$ is not in $C_i$ then $q$ is not in all $C_{j}$, $j>i$;
(3) $C_{i}$ and $C_{i+1}$ are side-gradually varied.
\end{definition}

We now show three reasonable definitions of simply connected spaces below.
We will provide a proof for the Jordan curve theorem under the third definition
of simply connected spaces. The Jordan theorem shows the relationship
among an object, its boundary, and its outside area.

Let $U_2$ be a subset of all minimal closed curves on $G=(V,E)$. Each element in
$U_2$ is a 2-cell.  $<G,U_2>=<V, E, U_2>$ defines a 2D topological structure discretely.

A general definition of a simply connected space should be :

\begin{definition} {\bf Simply Connected Surface Definition (a)} $<G,U_2>$ is
 simply connected if any two closed simple paths are homotopic.
\end{definition}

If we use this definition, then we may need an extremely long proof for
the Jordan curve theorem. The next one is the standard definition
which is the special case of the Definition 3. (Definition 3 is
too general, it is not needed here.)

\begin{definition} {\bf Simply Connected Surface Definition (b)}
A connected discrete space $<G,U_2>$ is  simply connected
if for any point $p\in S$, every simple cycle containing $p$ can contract to $p$.
\end{definition}

This definition of the simply connected set is based on the original meaning of
simple contraction. In order to
make the task of proving the Jordan theorem
simpler, we give the third strict definition of simply connected surfaces as
follows.

We know that a simple closed path (simple cycle) has at least
three vertices in a simple graph.
This is true for a discrete curve in a simply connected surface $S$.
For simplicity, we call an unclosed path an arc.
Assume $C$ is a simple cycle with clockwise orientation.
Let two distinct points $p, q\in C$. Let $C(p,q)$ be an arc of $C$ from $p$ to $q$ in
a clockwise direction, and $C(q,p)$ be the arc from $q$ to $p$ also in a clockwise
direction,
then we know $C=  C(p,q)\cup C(q,p)$. We use $C^{a}(p,q)$  to represent
the counter-clockwise arc from $p$ to $q$. Indeed, $C(p,q) = C^{a}(q,p)$.
We always assume that $C$ is in clockwise orientation.

\begin{definition} {\bf Simply Connected Surface Definition (c)}
A connected discrete space $<G,U_2>$ is simply connected
if for any simple cycle $C$
and  two points $p, q\in C$, there exists a sequence of simple cycle paths
$Q_{0},...,Q_{n}$ where $C(p,q)=Q_{0}$ and $C^{a}(p,q)=Q_{n}$ such that
$Q_{i}$ and $Q_{i+1}$ are side-gradually varied for all $i=0,\cdots, n-1.$.
\end{definition}

We have proved the following proposition in \cite{Che13}.
\begin{proposition}
       Definition (b) and Definition (c) are equivalent.
\end{proposition}

\section{The Jordan Curve Theorem in Discrete Space}

Since a simple cycle is a closed simple path that could be a surface-cell that cannot separate
$C$ into two disconnected components. So for the strict case of
Jordan curve theorem, we must use the closed discrete curve (not only a simple cycle).

A discrete curve $C$ does not contain a subset of vertices that are not all vertices of a 2 cell in $U_2$.
The intuitive meaning is that  $C$ does not contain any 2-cell.

A 2-cell will have two directions: the clockwise and the counter-clockwise. If we imagine a point at the center
of 2-cell, we will have two normals. (This is also true for 1-cell) This point will be called  the central pseudo point.
In the case of allowing the central pseudo points that will be called the Veblen point,
we will have the general Jordan Curve Theorem. We will prove this case
at the last of this section. The idea of  the central pseudo point  is valid when
embedding a surface (or cells)  into Euclidean space.

The proof of The Jordan Curve Theorem in Discrete Space will need to use the following proposition:

\begin{lemma}
       If $(x_0,x_1)$ is an edge, define $S(x_0,x_1)= S(x_0)\cup S(x_1)$. Then, $S(x_0,x_1)-\{x_0, x_1\}$ is a
simple path.
\end{lemma}

\noindent We can let $X=\{x_0, x_1\}$, and we denote $S(x_0,x_1)= S(X)$. However, this lemma is not true when
$x_0,x_1,\cdots,x_{k-1}$ is a discrete curve, $k> 2$. In general,

\begin{lemma}
 $S(X)-X$ is a degenerated (closed) simple path.
\end{lemma}

A degenerated (closed) path is a simple path with several unclosed discrete curve attach to some vertices on the path.
This is because, in discrete case, some part of the simple path shrink or folded into a poly-line.

When every ``angle'' in a discrete curve $C$ is little wide meaning that contains three 2-cells. In other words, for each three consecutive points
$x_{-1} x_0 x_1$ in $C$, a path from $x_{-1}$ to $x_{1}$ without passing $x_0$, there are must be two other vertices in between.
In addition, these two vertices are in a 2-cell or more 2-cells that does not contain $x_{-1}$ nor $x_{1}$. This means that
there is a 2-cell in between the neighborhoods of $x_{-1}$ and $x_{1}$.

We also request
any two vertices that are not adjacent will have the graph-distance greater than 2. The purpose is to guarantee that  $S(X)-X$ is
a simple closed path.

We denote the wideness of an angle is the minimum number of edges in the path from $x_{-1}$ to $x_{1}$ (each edge is in different 2-cell containing $x_0$).
Note that an angle with wideness 1 will make that $C$ is not a discrete curve. Since $C$ contains a triangle that is a 2-cell.)

\begin{lemma}
Let $C$ be a discrete curve and $X=\{x_0,\cdots,x_{k-1}\}$ be a path (arc) in $C$. ($x_0$ and $x_{k-1}$ are separated by at least three 2-cells.)
If the following conditions are satisfied:\newline
  (1) The wideness of each angle is 3 or greater, and \newline
  (2) For any two nonadjacent vertices $p$ and $q$ in $C$, any path not including edges in $C$ from $p$ to $q$ much contain three edges that belong to three different
    2-cells. \newline
\noindent then, $S(X)-X$ is a closed simple path.
\end{lemma}

We will prove this lemma in the proof of Theorem 1.

We know that we can easily make $C$ in 2D Euclidean plane to be a wide angle by adding some lines to the Veblen points on a edge. (Do not make extra
Veblen point used just for a 2-cell, use it for the 2-cell that share the Veblen point.)

The following theorem is for  $C$ with the wider angles.  To direct prove this theorem for general case meaning allowing that $S(X)-X$ is
a degenerated (closed) simple path.  We will give the proof in Appendix.

For a closed discrete curve, we have
\begin{theorem}({\bf The Jordan Curve Theorem in Discrete Space})
A discrete simply connected surfaces $S$
defined by Definition 5 (Definition (c)),
has the Jordan  property:
For a closed discrete curve $C$ on $S$,
if $C$ does not contain any point of $\partial S$,
$C$ divides $S$ into at least two disconnected components.
In other words, $S-C$ consists of at least two disconnected components.
\end{theorem}

\begin{proof}
Suppose that $C$ is a closed curve in a simply connected
surface $S$. $C$ does not reach the border of $S$, i.e.
$C\cap {\partial S} = \emptyset$.

We also check $C$ if it satisfies the condition of the wide angle for each triple consecutive vertices.
Also we add Veblen points and edges to make $C$ to be wide as necessary. So we can assume that:\newline
(1) The wideness of each angle is 3 or greater, and \newline
(2) For any two nonadjacent vertices $p$ and $q$ in $C$, any path not including edges in $C$ from $p$ to $q$ much contain three edges that belong to three different
    2-cells. \newline

Assume point $p\in C$, then suppose that $q$ and $r$
are two adjacent points of $p$ in $C$ with form of $...q p r,...$, where
the direction of  ... $q$ to $p$ to $r$ ...to $p$ is clockwise.
See Fig. 2
$\{p,r\}$ is a line-cell, then there are two 2-cells containing
$\{p,r\}$. Denote these 2-cells by $A$ and $B$ with clockwise orientation.

Our strategy is to prove that if
there is a point $a$ in $A$ which is not in $C$, and a point $b\in B$
and $b\notin C$, then any path from $a$ to $b$ must contain a point
in $C$. Then we can see that $S-C$ are not (point-) connected and we have the
Jordan curve theorem.

First, we want to prove that there must exist a point in $A-C$.
 If each point in $A$ is in $C$, since $A$ is a simple cycle,
then $C=A$. However, $C$ is not a surface-cell, so the statement can
not be true. Thus, there is a point $a\in A-C$. For the same reason
there is a point $b\in B-C$.   We assume
that $a$ is the last such point in $A$ starting with $p$, and $b$ is
the first such  point in $B$ starting with $p$. (see Fig. 2(a))
We always
assume  clockwise direction here for cell $A$ and $B$ unless we indicate otherwise.

\begin{figure}[h]
	\begin{center}

  \epsfxsize=4.5in
  \epsfbox{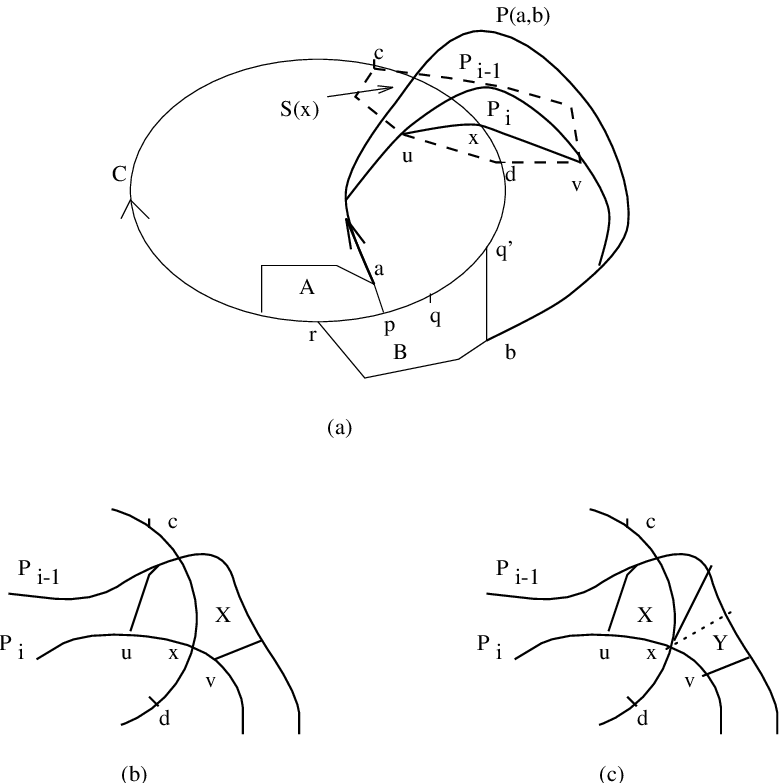}  

\caption{A close curve $C$ and the paths from $a$ to $b$}
\end{center}
\end{figure}

Let us make a summery of above idea: Suppose that $C$ is a closed curve.
(If $C$ is a closed simple path, we allow $C$ is a 2-cell and we
allow to assign a pseudo-point in the center of the 2-cell, we still can
prove this theorem.) The idea of the proof is to find two points in each
sides of the curve $C$.
This is because that for any 1-cell $(r,p)$ in $C$, there are two 2-cells $A$,$B$
sharing $(r,p)$ by the discrete surface definition.  $A$ must contain
a vertex $a$ and $B$ must contain $b$,
and they are not in $C$. $a$, $b$ are adjacent to some points in $C$,
respectively.
We are going to prove that from $a$ to $b$, any path must cross-over $C$.
That is the most important part of the Jordan curve theorem.

We assume, on the contrary, there is a simple path from $a$ to $b$
does not cross-over $C$, called
$P_(a,b)$ in Fig. 2 (a).
But we know there is $P(b,a)$ in $A\cup B$ (
i.e. $P(b,a)=b\cdots rp\cdots a$)  does cross-over $C$.
$P_(a,b)\cup P_(a,b)$ is a cycle in clockwise. (Fig. 2 (a) )

We know $S(r)$ is the neighborhood of $r$ in $S$.
So $S(r)$ contains all 2-cells containing
$r$. The boundary of $S(r)$ is a simple closed curve. (This is because
we always assume that $r$ is a regular point).
$a$ is on the boundary of $S(r)$.
(The boundary of $S(r)$ is $S(r)-\{r\}$). $A\cup B$ is a subset of $S(r)$.

We now prove that $P_(a,b)$ is not a subset of $S(r)$; otherwise, it must cross-over $C$.
(a 2-cell containing $r$ must
have an edge on $C$, or all points of the 2-cell are on the boundary of $S(r)$ except $r$).
If $P(a,b)$ does
not contain $r$, must be a part of boundary of $S(r)$ which is a cycle.
$r$ has two adjacent points on $C$, (If
they are not pseudo points, meaning here it can be eliminated or added on an edge
 that does not affect to the 2-cell) so these two points are also in the boundary
of  $S(r)$.
So there are only two paths from $a$ to $b$ on the boundary of $S(r)$.
These two points are not on the same side of the cross-over path containing $r$.
(The boundary of $S(r)$ was separated by
  the cross-over path containing $r$.) $P_(a,b)$ must contain such a point that is on $C$.

Therefore we proved $P_(a,b)$ is not a subset of $S(r)$.
Then $P_(a,b)\cup P(b,a)$ is a simple closed curve. ($P(b,a)$ passes $r$).
By the definition of the simply-connected surface,
there are
finite numbers of paths $P(a,b)=P_{0}(a,b)$,..., $P_{n-1}(a,b)$, such that
so that $P_{i}(a,b)$ and $P_{i+1}(a,b)$ are (side-)gradually varied.

In addition, $P_{n-1}(a,b)$ is gradually varied
to $P_{n}(a,b)=P^{a}(b,a)$ (a reversed $P(b,a)$ that passes $r$).

We now can assume that there is a smallest $i$ such that $P_{i}(a,b)$ cross over $C$,
but $P_{i-1}(a,b)$ does not.
(Fig. 2(a) ).   We will prove that
is impossible if $P_{i-1}(a,b)$ does not cross over $C$.

Let point $x$ in $P_{i}(a,b) \cap C$ and $x\notin P_{i-1}(a,b)$. There are two cases:
(1) cross over at a single point $x$ on $C$, or (2) cross over at
 a sequence of points on $C$. We will prove these two cases, respectively.


{\bf Case 1}: Suppose that $x=``p''=``q''$ in Definition 1 
(See Fig. 1 (a)(b)).
It means two curves $P_{i}(a,b)$ and $C$ share just one point $x$.
and assume $P_{i}(a,b)=...u x v...$ and $C= ... c x d ...$, where $v\ne d$.

We know that $u,v,c,d $ are in the boundary of $S(x)$, a simple cycle $S(x)-\{x\}$
(See Lemma 1).
There is a 2-cell $X$ (in between $P_{i-1}$ and $P_{i}$) contains $(u,x)$.
See Fig. 2 (b).

$X$ has a sequence of points $S1$ in $P_{i-1}$ and  a sequence of points $S2$ in $P_{i}$.
$X$ has at most two edges $e1$, $e2$ not in  $P_{i-1}\cup P_{i}$ ;
$S1$, $e1$, $S2$, $e2$,
are the boundary of $X$. $e1$ is the edge linking $S1$ to $S2$,
and $e2$ is the edge linking $S2$ to $S1$
counterclockwise.(Again, $e1$ may or may not be directly incident to $u$,
and $e1$ may be an empty edge if
$P_{i-1}$ intersects $P_{i}$ at point $u$. $e2$ may also in the same situation.)
We might as well assume that $x$ is the first
point on $P_{i}(a,b)$ (from $a$ to $b$ in path $P_{i}$ )that is in $C$.
Thus, $c,d\notin X$.
(If $c$ is in $X$ $c$ must be in $P_{i-1}$. if $d$ is in $X$, $x$ is not
only the cross over point. )

If $X$ contains $v$,  we will have a cycle $u\cdot d \cdot v (e2) (S1) (e1) $
in the boundary of $S(x)$
$(e2)(S1) (e1)$ contains only points in $P_{i-1}$ and $u$,$v$
(that are possible end points of $e1$, $e2$).
$c$ is on the boundary of $S(x)$ too. Where is $c$? It must be in the boundary
curves (of $S(x)$)
from $u$ to $d$ or the curve from $d$ to $v$.
Then $c,d$ in $S(x)$ must in the same side of $u x v$ which is part of $P_{i}$.
Therefore, $C$ and $P_{i}(a,b)$ do not cross-over each other at $x$.
See Fig. 2(b)
and the following extended figure.

If $X$ does not contain $v$, then there must be a 2-cell $Y$ (in between $P_{i-1}$ and $P_{i}$)
containing $(x,v)$. We can see that $X$ and $Y$ are line-connected in $S(x)$.
(See Fig. 2(c) and above extended figure.)
This is due to the definition of regular point of $x$, all surface-cells containing
$x$ are line-connected. Meaning there is a 2-cell paths they share a 1-cell
in adjacent pairs.

Since $X$ and $Y$ are line-connected, we can assume:

{\bf a)} $X$ and $Y$ share a 1-cell, i.e.  $X\cap Y= (x,y)$.  Then $y$
is on $P_{i-1}(a,b)$.
Let $e2$ be the possible edge from $v$ to $P_{i-1}(a,b)$. ($e2$ could be empty as $e1$)
and $u(e1)..y...(e2)v$ is on the boundary cycle of $S(x)$.  Except $u$ and $v$, $u(e1)..y...(e2)v$ is
on $P_{i-1}(a,b)$. $u...d...v$ is part of the boundary cycle of $S(x)$.
In addition, $c$ (that is not in $P_{i-1}(a,b)$) must be in the boundary curves (of $S(x)$)
from $u$ to $d$ or the curve from $d$ to $v$. Again,
$c,d$ in $S(x)$ must in the same side of $u x v$ which is part of $P_{i}$.
Therefore, $C$ and $P_{i}(a,b)$ do not cross-over each other at $x$.
(See Fig. 2(c) also see the above extended figure Fig. 3(c).)

\begin{figure}[h]
	\begin{center}

  \epsfxsize=5in
   \epsfbox{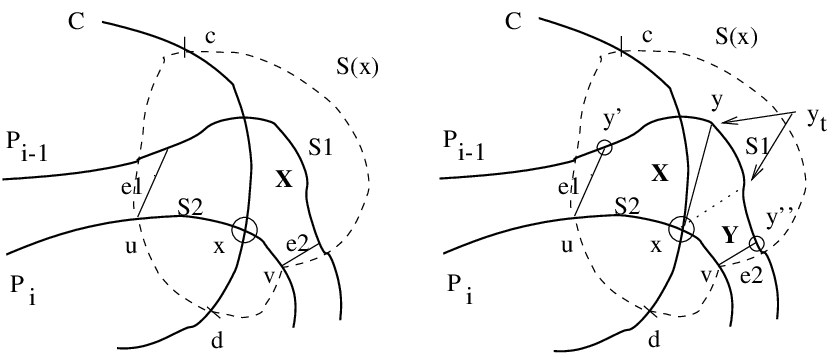 }
\caption{ Extended figures of Fig. 2 (b) and   Fig. 2 (c), respectively. }
\end{center}
\end{figure}

{\bf b)} $X$ and $Y$ share the point $x$, and there are line-connected 2-cells as a path in between $X$ and $Y$.
 $X\cap Y= x$ .  Let us assume that $e1$ incident to $P_{i-1}(a,b)$ at $y'$ ($y'$ is $u$ if
 $e1$ is empty ) and $e3$ incident to $P_{i-1}(a,b)$ at $y''$.
 We will have a set of points $y'=y_0, y_1,...,y_k=y''$ in $P_{i-1}(a,b)$. Each $y_t$  is contained
 in a 2-cell containing $x$. All $y_0, y_1,...,y_k$ are in the boundary cycle of $S(x)$.
 $c$ that is not in $P_{i-1}(a,b)$. $c$ must be in the boundary curves (of $S(x)$)
from $u$ to $d$ or the curve from $d$ to $v$. Thus,
$c,d$ in $S(x)$ must in the same side of $u x v$ which is part of $P_{i}$.
$C$ and $P_{i}(a,b)$ do not cross-over each other at $x$. (See Fig. 2 (c) and Fig. 3(c).)

{\bf Case 2:} Suppose $P_{i}(a,b)$ and $C$ cross over a sequence of points on $C$:
$P_{i}(a,b)=...u x_0 x_1...x_m v...$ and $C= ... c x_0 x_1...x_m d ...$, where $v\ne d$.

We still have $e1 = (u,y_0)$ and $e2=(v,y_k)$ where $y_0$ and $y_k$ are on
$P_{n-1}$ for some $k$.
Each $y_t$ , $t=0,1,...,k $, is  in a 2-cell that containing some $x_j$, $j=0,1,...,m$.

Note that : If $u$ does not have a direct edge linking to $P_{n-1}$, $u$ will be
in a 2-cell between $P_n$ and
$P_{n-1}$,  either $u$ is a pseudo point on $P_{n}$ for the deformation from $P_{n-1}$ to $P_n$, or
 $P_{n-1}$ and $P_n$ intersects at $u$. That $u$ is a pseudo point means here it has a neighbor that
 has an edge link to $P_{n-1}$, or the neighbor's neighbor, and so on. We can just assume here
 $u$ is the point that is adjacent to a point in $P_{n-1}$. In the theory, as long as
 $u$ is contained by a 2-cell such that all the points in the 2-cell are in $P_{n-1}$ or $P_n$.

\newpage

\begin{figure}[h]
	\begin{center}

  \includegraphics[width=3in] {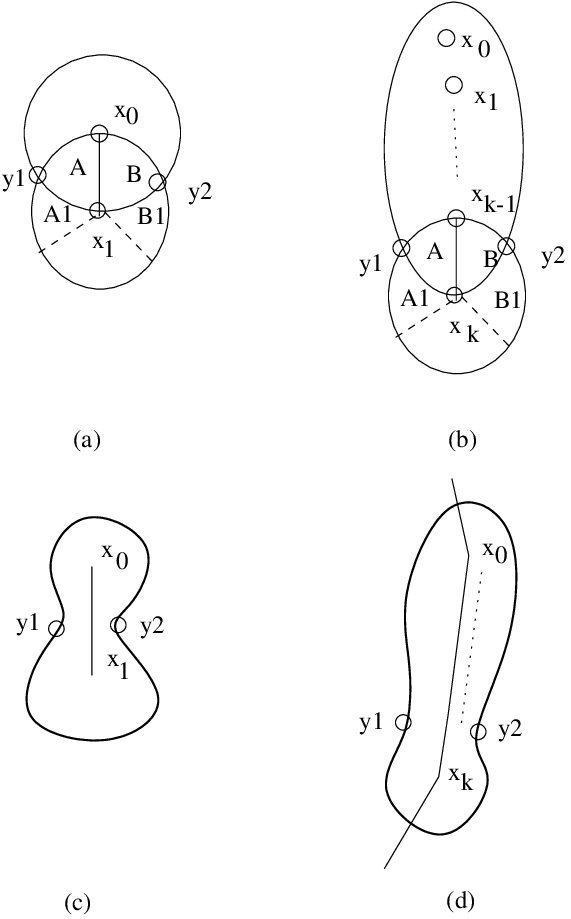}    

\caption {The union of neighborhoods of a sequence of adjacent points $S(x_0,...,x_k)$ and its boundary}
\end{center}
\end{figure}

The same way will apply to this case just treat $x_0$,...,$x_m$ to $x$ in Case 1.  We first get the union
of $S(x_0)$,...,$S(x_m)$. We want to prove that : The boundary of this union will be simple cycle too
under the condition of Lemma 5.

Using mathematical induction we can prove it. After that, we can prove the rest of theorem using the same method
presented in Case 1. See Fig. 4.  

The following is the detailed proof:   Let $S(x_0,...,x_k)=S(x_0)\cup...\cup S(x_k)$.

First, we will prove that the boundary of $S(x_0)\cup S(x_1)$ is a simple
cycle (it is a simple closed curve too). We know that $(x_0, x_1)$ is an edge in $C\cap P_{i}(a,b)$. Also, there are two 2-cells $A,B$ in $S(x_0)$ containing  $(x_0, x_1)$ .

$x_1$ is a boundary point in $S(x_0)$ , so no other 2-cell will contain $x_1$. In the same way, $S(x_1)$ also contains  $A,B$, and
$x_1$ is only contained in two 2-cells in $S(x_1)$. Therefore,  $S(x_0)\cap S(x_1)= A\cup B$ and  $A\cap B=(x_0,x_1)$.

Note that $A$ and $B$ are adjacent 2-cells. On the other hand, $x_1$ is on the boundary curve (that is closed) of $S(x_0)$.  So $x_1$ has two adjacent points
on this cycle, $y_1$ and $y_2$. (We assume that $y_1$ and $y_2$ are not pseudo points, so)  $y_1$ and $y_2$ are both on the boundary
of $S(x_0)\cup S(x_1)$. (If $y_1$ or $y_1$ is pseudo points, we can ignore $y_1$ or $y_2$ to find the a actual point that adjacent to $x_1$.)
$(x_1,y_1)$ has two 2-cells containing $(x_1,y_1)$ in $S(x_0)\cup S(x_1)$. For instance, in Fig. 4 (a) , $A$ and $A_1$ contain $(x_1,y_1)$ and
$ B$ and $B_1$ contain $(x_1,y_2)$. Thus, the boundary of $S(x_0)\cup S(x_1)$ is a closed curve that is formed by the arc from $y_1$ to $y_2$
 in the boundary of $S(x_0)$, plus  the arc from $y_2$ to $y_1$  in the boundary of $S(x_1)$.

Second, we assume the boundary of $S(x_0,...,x_{k-1})$ is a closed curve, when we consider the arc $x_0,...,x_{k-1}, x_{k}$ in $C$, we can prove
the boundary of $S(x_0,...,x_{k})$ is also a closed curve.

We know that we have two closed curves: Suppose that $Q$ is the boundary of $S(x_0,...,x_{k-1})$ , and $R$ is the boundary of $S(x_k)$. $(x_{k-1}, x_k)$ is in
$S(x_k)$, and  $(x_{k-1}, x_k)$ is in $S(x_0,...,x_{k-1})$ . There are two 2-cells $A$, $B$ containing  $(x_k-1, x_k)$ in $S(x_k)\cap S(x_0,...,x_{k-1})$.

$x_{k-1}$ is on the boundary cycle of $S(x_k)$, then $x_{k-1}$ must have two adjacent points in $R$, $y_1$, and $y_2$. $(x_{k-1},y_1)$ and $(x_{k-1},y_2)$ are two edges in $S(x_k)\cap S(x_0,...,x_{k-1})$. In the same way above, we will have the cycle passing $y_1$ and $y_2$ that is the boundary curve of $S(x_0,...,x_{k})$ .
This is because there is a 2-cell in between the neighborhoods of $x_{k-2}$ and $x_{k}$ as assumed by the wider angle on
each point on $C$,   $y_1$ is not a folding point, so are $y_2$. In other words, the point that enters $y_1$ in counterclockwise in $S(x_0,...,x_{k-1})-\{x_0,...,x_{k-1}\}$
differs from the point follows $y_1$ in counterclockwise in $S(x_k)-\{x_{k}\}$.

Thus, we have proved that  the boundary curve of $S(x_0,...,x_{k})$  is a simple closed curve.  In the rest of the proof, we will treat $S(x_0,...,x_m)$ to be $S(x)$ in Case 1.
See Fig. 5.
We now use denote  $X=\{x_0,...,x_m\}$  and $X$ is an arc in $C$. (Please note that in Case 1, $X$ was used as a 2-cell. Now $X$ is an arc in $C$.  )

\begin{figure}[h]
	\begin{center}

    \includegraphics[width=3.5in] {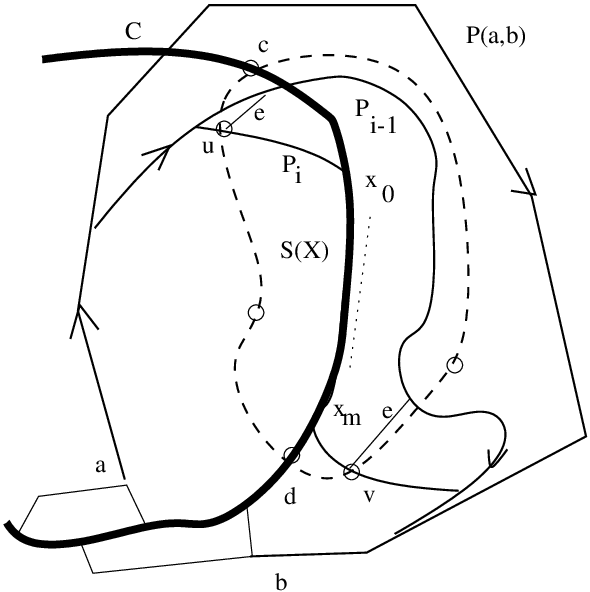}  

\caption{A pair of crossing-over curves $P_{i}$ and $C$  pass an arc $X=\{x_0,...,x_m\}$ }
\end{center}

\end{figure}

In the rest of the proof,  we will prove:
if $P_{i-1}(a,b)$ and $C$ are not cross over each other, then, $P_{i}(a,b)$ and $C$ will not be cross over each other.
Therefore, any $P(a,b)$ must cross over $C$. This completes the proof of the discrete Jordan curve theorem.

Let us first state again that $P_{i}(a,b)$ passes $x_0...x_m$ but $P_{i-1}(a,b)$ does not contain any point of $\{x_0,...,x_m\}$.
In addition, $P_{i-1}(a,b)$  and $P_{i}(a,b)$ is gradually varied, i.e. $P_{i}(a,b)$ was deformed from $P_{i-1}(a,b)$ directly. We also know that
$S(X)=S(x_0,...,x_m)$ is the neighborhood of the arc in $C$, i.e. the arc $x_0,...,x_m$ is a part of the closed curve $C$. The boundary of
$S(X)=S(x_0,...,x_m)$ is a closed curve too.

$u, v,  c,  d$ are on the boundary of  $S(x_0,...,x_m)$  (Assume $u, v,  c,  d$ are not pseudo points, otherwise, we can find
corresponding none-pseudo on  the boundary of  $S(x_0,...,x_m)$.)  $u, (x_0,...,x_m),v$  is a part of $P_i$
We also know that $c$ and $(x_0,...,x_m)$ are not in $P_{i-1}$. There will be two 2-cells,
$U$ and $V$, are in between $P_{i}(a,b)$ and $P_{i-1}(a,b)$ (all points of $U$ and $V$ are in $P_{i}(a,b)\cup P_{i-1}(a,b)$)  such that $(u, x_0) \in U$ and $(x_m,v)\in V$.

Let $P_{i-1}\cap U = S1$ and $P_{i}\cap U = S2$. Let $e1$ be the edge in $U$ linking $S1$ to $S2$ (in most cases, $e1$ incident to $u$, but not necessarily ), and let $e2$ be the edge in $U$ linking $S2$ to $S1$ (possibly starting at $x_0$).
So, $(e2) (S1) (e1) (S2)$ are the boundary of $U$, counterclockwise.

{\bf Subcase (i):}  If $U$ contains $v$ ($U=V$), all points in $U$'s boundary are contained in $S(\{x_0,...,x_m\})$ by the definition of $S(x_0)$.  we will have a cycle $u\cdot d \cdot v (e2) (S1) (e1) $ in the boundary of $S(X=\{x_0,...,x_m\})$.
 $c$ is on the boundary of $S(X)$ too. But $c\notin P_{i-1}$.  It must be in the boundary curves (of $S(X)$)
from $u$ to $d$ or the curve from $d$ to $v$.
Then $c,d$ in $S(X)$ must in the same side of $u X v$ which is part of $P_{i}$.
Therefore, $C$ and $P_{i}(a,b)$ do not cross-over each other at $X$. (See Fig. 5.)

{\bf Subcase (ii):} If $U$ does not contain $v$, then there must be a 2-cell $V$ (in between $P_{i-1}$ and $P_{i}$) containing $(x_m,v)$.

Let $e1=(p1,p2)$ be the edge in $U$ incident to a point in $P_{i-1}$ and a point in $P_i$, respectively.
(In most cases, $e1$ incident to $u$, i.e. $u=p2$, but not necessarily ). And let $e2=(r2,r1)$ be the edge in $V$ incident to a point in $P_{i}$ and a point in $P_{i-1}$, respectively. $r2$ is usually $v$.

Note that: $c$ must not be in $U$. Gradual variation
(direct deformation) means that  each point in each 2-cell of $U$ and $V$ in between $P_i$ and $P_{i-1}$  must be in $P_i \cup P_{i-1}$. Formally,
$P_i$ $XoRSum$  $P_{i-1}$ is a set of 2-cells; every point in these 2-cells is in $P_i \cup P_{i-1}$ .

We can see that $U$ and $V$ are line-connected in $S(X)$ by the definition of line-connected paths,  meaning there is a path of 2-cells where each adjacent pair shares a 1-cell.
(See Fig. 4(d)
and Fig. 5.)

From $r1$ to $p1$, there is an arc in $P_{i-1}$ . To prove that all points in this arc are in the boundary of $S(X)$ we need
to prove each point on the arc must be in a 2-cell that contains a point in $\{ x_0,..., x_m\}$, and this 2-cell is
other than (except this 2-cell is) $U$ or $V$. It gives us some difficult to prove it.  The above discussion seems not very productive.

We found a more elegant way to prove  this case by finding another simple path (or pseudo curve) that cross-over $C$.  The method is the following:
If $U\neq V$,  there must be a $x_k$ in $\{ x_0,..., x_m\}$, $x_k$ has an edge linking to $P_{i-1}$. (Otherwise, $u,  x_0,..., x_m, v$
are in a 2-cell that contains some points in  $P_{i-1}$. Therefore, $U=V$.) We can also assume that $k$ is not $m$, otherwise,
$v$ is in $P_{i-1}$, so $U=V$.  See Fig. 6.

\begin{figure}[h]
	\begin{center}

  \epsfxsize=3.5 in
   \epsfbox{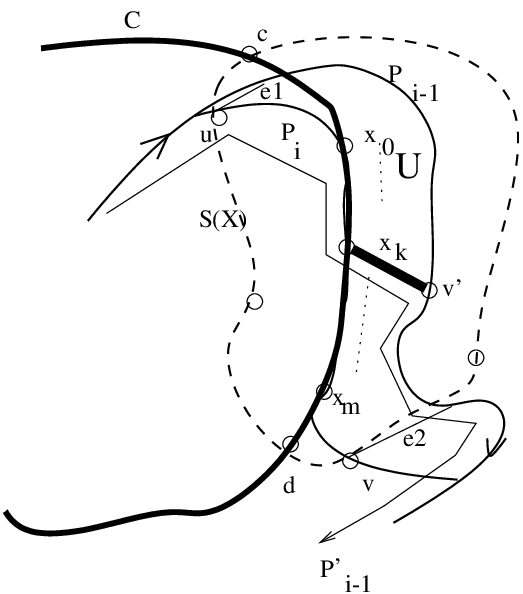} 

\end{center}
\caption{ An edge only starts at $x_k$ linking to $P_{i-1}$;  $x_0,...,x_{k-1}$ do not have any edge to $P_{i-1}$. }

\end{figure}

We select the smallest $k$ having an edge linking $x_k$ to $P_{i-1}$, $0\le k\le m-1$.  $x_k$ is in both $P_{i}$ and $C$.
We might as well let $(x_k,v')$ is such an edge, and $v'$ is a point in $P_{i-1}$. Therefore, we will have the new path (simple path), $P'_{i}$.
This new path has two parts: The first part  is the same as   $P_{i}$
before and including the point $x_k$, and the second part is the partial path (curve) of $P_{i-1}$ after point $v'$ . This path $P'_{i}=...,u, x_0, ...,x_k, v',...$ does cross-over $C=...,c, x_0, ...,x_k,x_{k+1},...,x_{m}, d,... $ .
It is obvious that $P_{i-1}, P'_{i}, P_{i}$ are gradually varied.
This is because we just inserted a path in between of $P_{i-1}$ and $P_{i}$.

This new path $P'_{i}$ has such a good property that is
$\{x_0, ...,x_{k-1}\}$ do not incident (link) an edge that has an end vertex in $P_{i-1}$ .
Since $v'$ is in $P'_{i}$,  the 2-cell $V'$ (in $P_{i-1}$ $XoRSum$ $P'_{i}$ {\footnote {Here $XoRSum$ is exclusive ``OR'' operation defined in \cite{Che04}, just like $sum(modulo2)$ in \cite{New}}} )
contains $v'$ also contains $(x_{k-1}, x_k)$ and $(x_k, v')$ (in $S(x_k)$).
Since no edge from $x_0,...,x_{k-1}$ to $P_{i-1}$ , $U$ containing $u$ is just $V'$.
We will have just {\bf Subcase (i)} using $P'_{i}$ to replace $P_{i}$.

The entire theorem is proven.   \qed
\end{proof}

\section{The Jordan Curve Theorem for Generalized Simple Closed Paths}

The discrete Jordan curve theorem proved in last section has a little difference from the classical description of The Jordan curve Theorem.
This is because that discrete curve has its own strict property: $C$ does not contain any 2-cell. In order to satisfy the classical
form. We need to use central pseudo points, we call it the Veblen point, for each type of cells, especially 1-cells (line-cells) and 2-cells (surface-cells)
So we will allow the simple path (semi-curve) in the proof of  the Jordan curve theorem. In fact, a little modification will assist the proving of the theorem.
The rest of work is just to  prove that there are only two (connected) components in $S-C$.

A 2-cell in this paper is a small unite, the smallest unite that does not contain any other 2-cell.
A 2-cell in discrete space contains a central pseudo point that is called the Veblen point in this paper. At this point,
we can define two normals, one is in clockwise direction and another is is on counterclockwise direction. We can realize
it in Euclidean plane.

\begin{theorem} ({\bf The Jordan Curve Theorem for Generalized Simple Closed Paths}) Let $S$ be a discrete simply connected surfaces, ($S$ can be closed or a discrete plane embedded
in 2D Euclidean Space). A closed simple path (0-cell connected semi-curve) $C$ which does not contain any point of $\partial S$
divides $S$ into two components (in terms of allowing central pseudo points for each cell). In other words, $S-C$ consists of two components. These two components are disconnected.
\end{theorem}

\begin{proof}
In this proof, we can put the central pseudo points (the Veblen points) for each  1-cells and 2-cells to assist our proof. The reason is that if we
embed  1-cells and 2-cells into Euclidean plane or higher dimensional space. We can always find the  central points for each
cell.   The idea of the central pseudo points is at least valid in Euclidean space. In fact, the central pseudo points also have two normal directions
for a 2-cell. It also has two directions for a 1-cell. For instance, $\{a,b\}$ is an edge, $a\to b$ and $b \to a$ are two directions.

In the proof of Theorem 1,  we know that we have two 2-cells $A$ and $B$ at the different side of  the cycle $C$.
We proved that point $a\in A$ and $b\in B$ are not connected in $S-C$.

$C$ has the orientation of  clockwise (or counterclockwise as we first made). $(p,r)$ is  clockwise
in $A$, but $(p,r)$ is counterclockwise in $B$. So we call $A$ is  clockwise, and $B$ is  counterclockwise.  For each edge $e_i$ (e.g. $(p,r)$) in $C$,
we will have two 2-cells  containing $e_i$, denoted by $A_i$ and $B_i$. There must be one in  clockwise and another is in counterclockwise.

We always assume that $A_i$ is  clockwise and  $B_i$ is counterclockwise. We now add all the central pseudo points to all 1-cells and 2-cells in $S$.
And immediately remove all central pseudo points from 1-cells in $C$. (This operation is to stop  a path will go through the  central pseudo points
on $C$.)

We also know in our assumption: each 2-cell must have at least three edges (1-cells) in its boundary. This is because $S$ is a simple graph. (We
can always add a point to make it in Euclidean space.)  We also assume that $C$ does not reach the border of $S$ meaning that $C\cap {\partial S}=\emptyset$.

{\bf Case 1:  A special case}  $C$ is the boundary of a single 2-cell, denoted as $A$.    A simple path could be just the boundary of
a 2-cell. In this case, we have a central pseudo point in the cell $A$.  So this theorem is virtually true if we can prove that  $S-C$ is point-connected. That is to say that except
the  central pseudo point in $A$, $S-C$ is a point-connected component. In other word, $S-A$ is one point-connected component.

We can prove this  is because of the following facts:  For each cell $B$ in $S-A$, if $B$  has an edge in $C$. The boundary of  $B$
is a simple closed path. If the boundary of $B$ is $C$, then $S=A\cup B$ since every edge has shared by two 2-cells ($A$ and $B$) already.
We have two components in $S-C$: the central pseudo point of $A$ and  the central pseudo point of $B$.

{\bf Case 2:  The general case}

We know 2-cells that are joint with an edge in $C$ have two types: the clockwise type, denoted as $A_i$, $i=1,\cdots, m$, $m> 0$, and the counterclockwise type, denoted as $B_j$, $i=1,\cdots, n$, $n> 0$.

We can prove that all $A_i$ are connected without using points in $C$. This is because
that any point $p$ in $C$ is contained by two 1-cells $e1$ and $e2$ in $C$. These two 1-cells are contained by $A_i$ and $A_j$,
respectively. If $A_i =A_j$, it is connected.   If A If $A_i$ and $A_j$ share an edge, then, the central pseudo points of $A_i$ and $A_j$ are connected.
If  $A_i$ and $A_j$ do not share an edge, we know $A_i$ and $A_j$ are in $S(p)$, there must be a cycle contains some edges in $A_i$ and
some edges of $A_j$, and $e1 \cup e2$. So $A_i$ and $A_j$ are connected (meaning through their central pseudo points ) do not pass $e1 \cup e2$.
Therefore, all $A_i$'s (meaning using their central pseudo points) are connected. In other words,  $e1 \cup e2$ split $S(p)$ into two parts, one called $Part_A(p)$ include
$A_i$ and $A_j$, and another one, $Part_B(p)$ include some $B$'s. $A_i$ and $A_j$ are point connected in  $Part_A(p)$ without passing any point in  $e1 \cup e2$.
All cells that are not $A_i$ or $A_j$ in $Part_A(p)$ will also assign as the  clockwise type, i.e. $A_k$ for some $k$. So all $A_i$'s are connected.

In the same way, we can prove that all $B_i$'s are connected. (  $C\cap {\partial S}=\emptyset$.)

We now prove that any point $x$ in $S-C$, must be connected to  the component containing $A_i$ or to the component containing $B_i$.
We know that any two points are point-connected by a path in $S$. Let $c\in C$, $P(x,c)$ is such a path connecting $x$ and $c$.
Note that every point $c'$ in $C$ is contained by $S(c')=Part_A(c')\cup Part_B(c')$.  Since $x$ is not in $C$,  $P(x,c)$ (which has
finite numbers of points) must contain the first point in $C$, we assume it is $c'$. In many cases, $c'=c$. Let $P(x,c)= x\cdots x' c' \cdots c$.
Then $x'$ must be not in $C$. Thus, $x'$ must be in $S(c')-C$. $x'$  must be in some $A_i$ or$B_j$ because  $S(c')=Part_A(c')\cup Part_B(c')$.

In other words, there must be a first point
in $P(x,c)$, $x'$, that is adjacent a point $c'\in C$ ( may or may not be point $c$). $(x', c')$ must belong to an $A_i$ or $B_j$. So
if $(x', c')$ belong to $A_i$, $x$ is a point connected to the central pseudo points of $A_i$ . We call it component ${\cal A}$.  All points in ${\cal A}$
are connected since $A_i$ are connected for all $i$.

If $(x', c')$ belong to $B_j$, $x$ is a point connected to the central pseudo points of $B_j$ . We call it component ${\cal B}$. All points in ${\cal B}$
are connected since $B_j$ are connected for all $j$. We also know $S-C ={\cal A}\cup {\cal B}$ since $x$ was selected from  $S-C$.

According to the proof of on Theorem 1, there is $a$ in some ${A_i}-C$ is not connected to $b$ in some $B_j-C$ in $S-C$.
 $a$ is in  ${\cal A}$, and $b$  in ${\cal B}$.  Therefore, any point in ${\cal A}$ is not connected to any point in ${\cal B}$ in
in $S-C$.  (Otherwise, $a$ will be connected to that point in  ${\cal A}$, then will be connected to any point in  ${\cal B}$.)

 We now complete the proof of Theorem 2,
 the general Jordan Curve Theorem. \qed
\end{proof}

\section{Subdivision of Triangles and Jordan Curve Theorem in Euclidean Space}

In Theorem 2,  we allow the simple path (pseudo-curves) for the Jordan Curve Theorem.
This is the general case of Jordan Curve Theorem in discrete space. We know that
2D Euclidean space can be partitioned into triangles and it is simply connected in this discrete space in terms of simplicial complexes.
So we can prove  the Jordan Curve Theorem for 2D Euclidean plane.

The only problem is that we need to assume that
we have a refinement process that will make the triangulation (joint with the simple curve $C$) infinitively approximates $C$.

The following construction will bring a more satisfied answer. We will use the mid point subdivision method
to refine a triangulation. Even though we can use barycentric subdivision to refine a triangle, but it is computationally expensive.
In the following figure, the mid point subdivision method will partition the triangle into four small triangles if we
agree a curve can be represented as $C(t)$, $t\in [0,1]$. then we can get $C(0.5 \cdot t)$, And so on so fourth.

\begin{figure}[th]
	\begin{center}

  \includegraphics[width=3.5in] {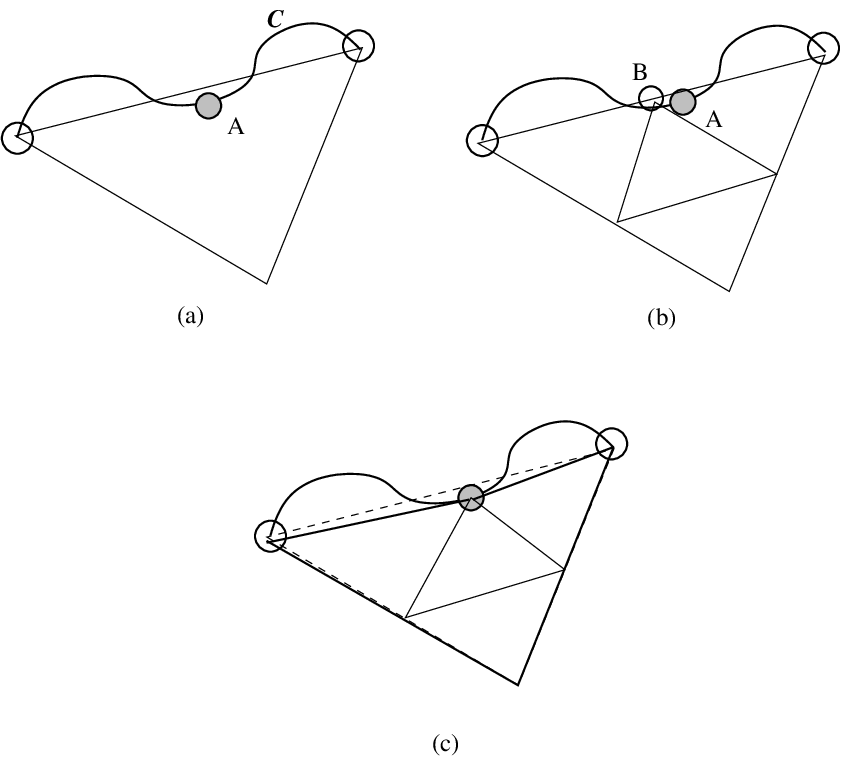}

\end{center}
\caption{The mid point subdivision method to partition a triangle:(a) The curve $C$ and its approximation triangle
         (b) Subdivision of the triangle, but the mid point of the edge of triangle is not
         exactly on the mid point of the arc of $C$, (c) make the two mid points as one on the curve $C$. }

\end{figure}

{\bf Algorithm 1}  The mid point subdivision method to partition a triangle along with the boundary curve $C$.

  Step 1: The curve $C$ joints with a triangle at two points.  (Fig. 7 (a))

  Step 2: Subdivision of the triangle using the mid point method. The mid point of the edge joining with $C$ in triangle is not
         exactly on the mid point of the arc of $C$.   (Fig. 7 (b))

  Step 3: Make the two mid points as one on the curve $C$. And modify the subdivision triangle accordingly. (Fig. 7 (c))

(This algorithm might not be a new algorithm. It is a natural way to do it. If someone already found this algorithm,
We will cite his/her work.)

If the close (boundary) curve $C$ is a continuous curve, the process of making the mid point subdivision is always valid.
This process is only for the approximating that will cover all points. This theorem is valid for the point of approximation.
For any $\epsilon$, we can find a triangulation where the length of each edge is smaller than $\epsilon$, and the
vertices on the boundary (in discrete term $B_C$ are the points on the curve $C$.

There is no way to link outside of the $B_C$ without passing $B_C$.  $B_C$ separates the plane into two
components. $B_C$ is infinitively close to $C$.  So the theorem  is proven in the continuous case.

The second method could be the following:
Can we define a line with with, then we can make the line thinner and thinner to approximate the answer? What we can
state is that for any $C$ and $\epsilon$, we can find a $B(C,\epsilon)$ that is the approximation of $C$ with respect to
$\epsilon$. Each edge is shorter than $\epsilon$.

(A mapping from $B(C,\epsilon)$  to $C$, and the distance between two adjacent points on $C$ is bounded by
a bounded function of $\epsilon$ since $C$ is continuous)

So we can prove that every curve with width $\epsilon$, from inside of curve $C$ to outside of curve $C$, must
pass a vertex in $B(C,\epsilon)$. We already proved this.  When  $\epsilon$ goes infinitively small, we have
this theorem proved. This is one way of thinking.

The Jordan Curve Theorem might only have some ``approximation'' proofs
along with the progress of mathematics. However, for a safe part, finite and discrete proofs of this great theorem is
very important.\\

{\bf Acknowledgment} The author would like to express many thanks to Professors Feng Luo and Xiaojun Huang at Rutgers University.
They have provided many helpful comments. Dr. Luo mentioned a result on triangled surfaces: It stats that every two poly-line are deformable to
each other by using a sequence of poly-lines where each pair of adjacent two only differs by a triangle. If this result is directly used
the proof of Theorem 1 in this paper, the proof would be much easier for a triangulated surface.   Due to the nature of the difficulty of the
proofs on the Jordan curve theorem historically, the author will hold total responsibility of this paper when an error occurs in the paper. \\

\section {Appendix  A: The Direct Proof of JCT in Discrete Cases}

The following is the proof of Theorem 1 without adding any new 2-cells for avoiding a degenerate simple path. In other words, we
allow that $S(X)-X$ is a degenerate simple path in the direct proof.

We have add some 2-cells in order to preserve the boundary cycle of $S(X)$ in Lemma 5 and Theorem 1. However, the Jordan curve theorem is
true for the any closed discrete curve $C$ on surface $S$ defined in this paper.

Recall the boundary of  $S(X)$,  $S(X)-X$. We can say that even though $S(X)-X$ might not be a closed simple path. But it is a degenerate simple
path (Fig. 8)

\begin{figure}[th]
	\begin{center}

  \includegraphics[width=3in] {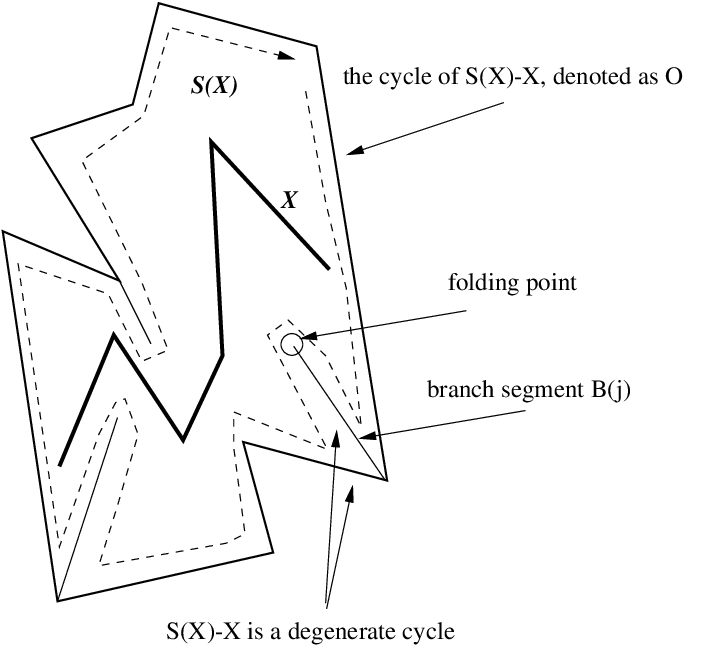}

\end{center}
\caption{$S(X)-X$  is a degenerate simple path in $S(X)$. $S(X)-X$ consists of two parts: the boundary cycle $O$ and several branch line segments $B_j$.}

\end{figure}

Our purpose is to prove that the point $c$ in Fig.2(a) is on the cycle $O$ not on the branch $B_j$ of $S(X)-X$ , or there is a point on $O$ that has the similar functionality like $c$ . Therefore, we can still make the same statement as it in the proof of
Theorem 1.

In the proof of Lemma 5, we can see for general case, $y_1$ might be a folding point since it is possible that $S(x_k)$'s boundary attached with
the boundary of $S(x_0,\cdots,x_{k-1})$.   See Fig. 4 (b) and (d). We let $X_{i}=\{x_0,\cdots,x_{i}\}$ for convenience later.

There are two cases we are interested in: (1) $c$ is the last point of branch point $B_j$ respect to $O$, (2) $c$ is the middle point of branch point $B_j$. If $c$ is on the cycle (the first point of $B_j$ to $O$ is also in $O$), we will use this case in the current proof which will be similar to the proof in Theorem 1.

{\bf Case 1}: If $c$ is the last point of branch point $B_j$, we know that  the last point of branch $B_j$ （which is $c$） is adjacent to a $x_i$ in $X$,  $i\ne 0$ . since $c$ is
also adjacent to $x_0$. So $\{c,x_0,\cdots,x_i\}$ is already a cycle. So $C=\{c\}\cup X$. This is not possible since we only consider $X$ is at most having $|C|-2$ points
when considering cell $A$ and $B$, there is an edge on $C$ not in $X$ . See Fig. 2. (a).

{\bf Case 2}: If $c$ is in the middle of $B_j$, there is an edge containing $c$ in $B_{j}$ toward to $b_0$ on $O$, the first point on $B_{j}$ . There
are two 2-cells ($E$ and$ F$) containing this edge in $S(X)$. On the other hand, $c$ has two neighbors in $C$: one is $x_0$ another denotes $x_{-2}$ . Consider $S(c)$,
$E$ and $F$ are also in $S(c)$. See Fig. 9. , we only have four possible subcases. $E$ is not reachable for $X$ in Fig  (b) and (c),   and $F$ is not
reachable from $X$ in Fig. 9 (a) and (d).
By the definition of $E$ and $F$, they are 2-cells that contain a point in $X$. Therefore, $x_{-2}$ must be on  $B_{j}$. If $x_{-2}$ is not $b_0$, we can do it again to
find $x_{-3}$ in $C$ that is also on $B_{j}$. So on so forth,  we will have one point $x_{-k}$ on $C$ that is $b_0$. $b_0$ is on $O$. Please note: $c=x_{-1},x_{-2},\cdots,x_{-k}$ are in $C$.

\begin{figure}[th]
	\begin{center}

  \includegraphics[width=4in] {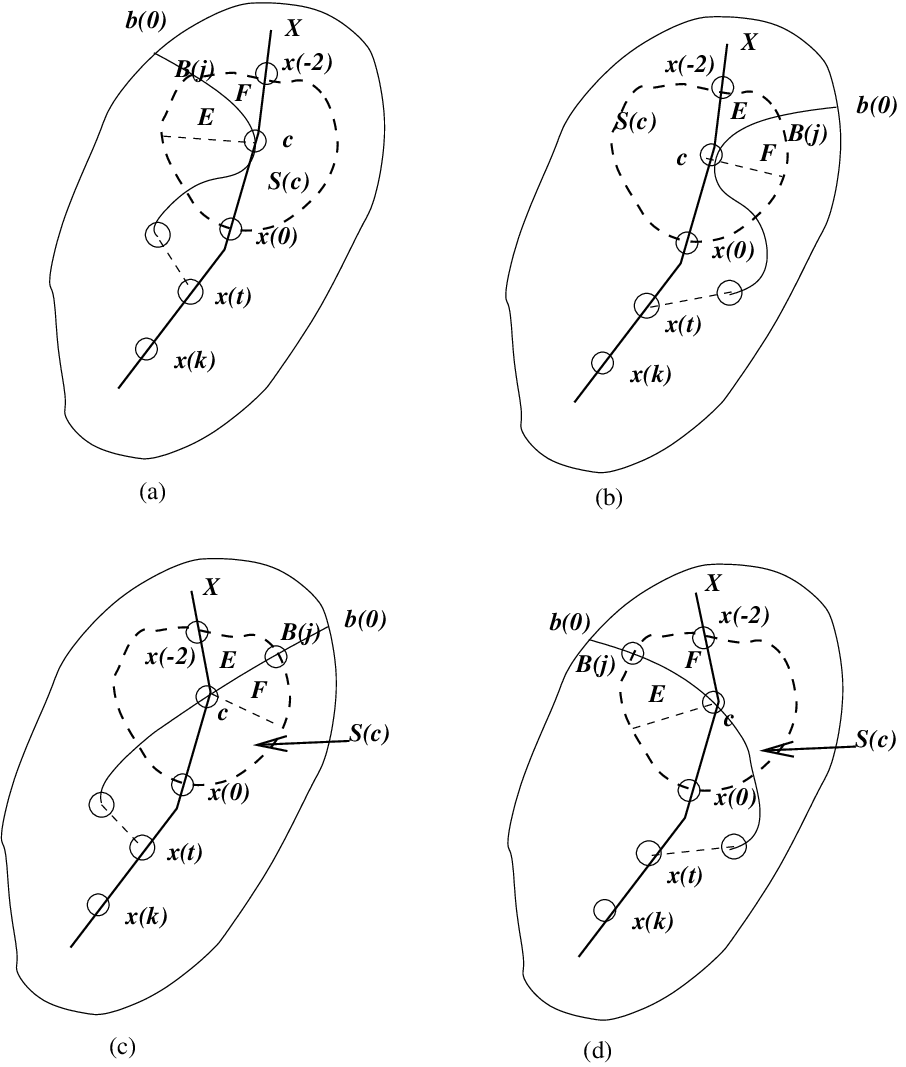}

\end{center}
\caption{ Only four cases when $c$ is in the middle of $B_{j}$.}

\end{figure}

We have proved that there is a point $c_0$ (it was denoted as $x_{-k}$) near $x_0$, $c_0$ is on $O$.  ($c_0$ is $c$ in the most time). The same thing will be true for $d$, denote the points as $d_0$.
In $S(X)-X$, from $c$ to $c_0$, every point is in $C$ as we proved above.   Therefore the situation is exactly the same as we discussed in the proof of Theorem 1.

Up to now we should be able to declare that we have proved completely Theorem 1 without adding any extra 2-cells in the proof of Theorem 1.

In order to walk through the whole proof with each small detail. We now repeat the proof of the case $P_{i}(a,b)\cap C =X$ where $X=\{x_0,...,x_m\}$ ($m>0$)  as in Theorem 1.

In the rest of the proof,  we will prove: If $P_{i-1}(a,b)$ and $C$ are not cross over each other, then, $P_{i}(a,b)$ and $C$ will not be cross over each other.
Therefore, any $P(a,b)$ must cross over $C$. This completes the proof of the discrete Jordan curve theorem.

Let us first state again that $P_{i}(a,b)$ passes $x_0...x_m$ but $P_{i-1}(a,b)$ does not contain any point of $\{x_0,...,x_m\}$.
In addition, $P_{i-1}(a,b)$  and $P_{i}(a,b)$ is gradually varied, i.e. $P_{i}(a,b)$ was deformed from $P_{i-1}(a,b)$ directly. We also know that
$S(X)=S(x_0,...,x_m)$ is the neighborhood of the arc in $C$, i.e. the arc $x_0,...,x_m$ is a part of the closed curve $C$. The boundary of
$S(X)=S(x_0,...,x_m)$ is a degenerate simple path shown in Fig. 8.

$u, v,  c_{0} [c],  d_{0} [d] $ are on the boundary of  $S(x_0,...,x_m)$  (Assume $u, v,  c_0,  d_{0}$ are not pseudo points meaning only have two neighboring points in $S$, otherwise, we can find
corresponding none-pseudo on  the boundary of  $S(x_0,...,x_m)$.)  $u, (x_0,...,x_m),v$  is a part of $P_i$
We also know that $c_0 [c]$ and $(x_0,...,x_m)$ are not in $P_{i-1}$. There will be two 2-cells,
$U$ and $V$, are in between $P_{i}(a,b)$ and $P_{i-1}(a,b)$ (all points of $U$ and $V$ are in $P_{i}(a,b)\cup P_{i-1}(a,b)$)  such that $(u, x_0) \in U$ and $(x_m,v)\in V$.

Let $P_{i-1}\cap U = S1$ and $P_{i}\cap U = S2$. Let $e1$ be the edge in $U$ linking $S1$ to $S2$ (in most cases, $e1$ incident to $u$, but not necessarily ), and let $e2$ be the edge in $U$ linking $S2$ to $S1$ (possibly starting at $x_0$).
So, $(e2) (S1) (e1) (S2)$ are the boundary of $U$, counterclockwise.

{\bf Subcase (i):}  If $U$ contains $v$ ($U=V$), all points in $U$'s boundary are contained in $S(\{x_0,...,x_m\})$ by the definition of $S(x_0)$ (since $(u, x_0) \in U$ ).  we will have a cycle $u\cdots d_0 [d] \cdots v (e2) (S1) (e1) $ in the boundary of $S(X=\{x_0,...,x_m\})$.

$c_0$ and $c$ are on the boundary of $S(X)$ too. Especially,  $c_0$ is on $O$ in $S(X)-X$(See Fig. 8). In addition, $arc(c,c_0)$ in $C$ is in $S(X)-X$.  But $c_{0}\notin P_{i-1}$.  $c_0$ (and $arc(c,c_0)$)  must be in the boundary path of $S(X)$
from $u$ to $d_0$ or from $d_0$ to $v$. $c_0$ is on $O$. So if $c_0$ is in between $u$ to $d_0$, see Fig. 10,
then $c_0,c,X,d,d_0$ (an arc in $C$) in $S(X)$ must be in the same side of $u X v$ which is part of $P_{i}$. The same reason will
apply to the case that $c_0$ is in between $d_0$ to $v$. Please also note that $u$ and $v$ are in the boundary (degenerate) path as shown in Fig. 8 not necessary on $O$.
Therefore, $C$ and $P_{i}(a,b)$ do not cross-over each other at $X$. (See Fig. 10.)

\begin{figure}[th]
	\begin{center}

  \includegraphics[width=3in] {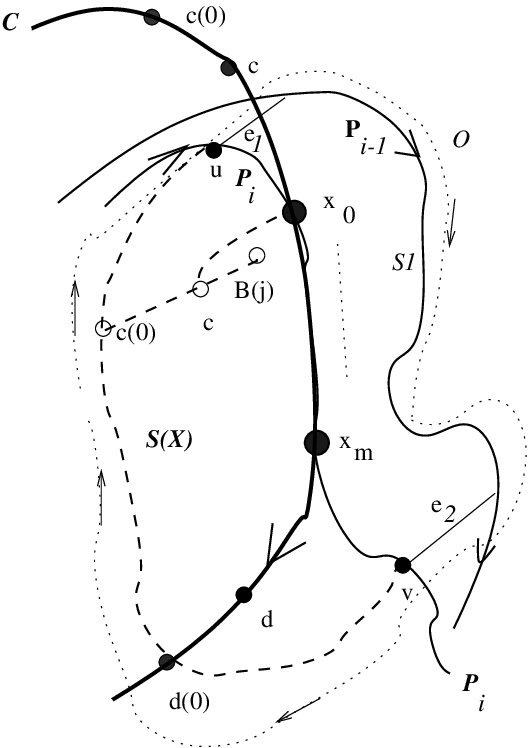}

\end{center}
\caption{The more general case as shown in Fig 5 where $S(X)-X$ is a degenerate simple path.}

\end{figure}

{\bf Subcase (ii):} If $U$ does not contain $v$, then there must be a 2-cell $V$ (in between $P_{i-1}$ and $P_{i}$) containing $(x_m,v)$.
We still want to find another simple path (or pseudo curve) $P'_{i}$ that cross-over $C$ in between $P_{i-1}$ and $P_{i}$.
And between $P_{i-1}$ and $P'_{i}$, there is only 2-cell just like {\bf  Subcase (i)}. The construction method is exactly the same as
{\bf Subcase (ii)} in the proof of Theorem 1 in Section 2.  Therefore, we can obtain the case is just like {\bf  Subcase (i)} above.
The proof is completed.


\section {Appendix  B: Simply Connected Space and the Euclidean Plane}

In this appendix, we will prove that the Euclidean plane is simply connected under our definition. As a by-product, this proof will also
prove  the Jordan–Schoenflies theorem.

The Jordan–Schoenflies theorem states that: A simple closed curve on the Euclidean plane separate the plane into two connected components. The component with bounded closure
is homeomorphic to the disk.

For a given simple polygon $C$ in the plane $R\times R$, our proof is to make a standard triangulation (the equilateral triangle partition) on the plane. That will dense enough such that
there is no two and more vertices of $C$ are inside of a triangle or an edge.  Then we modify the triangulation so that $C$ will be the boundary of the triangulation.
Then, we design an algorithm that will do contraction of $C$ to any given point on $C$. This contraction is a sequence of side-gradually varied paths.
So we can prove   $R\times R$  is simply connected under our definition of this paper.

Since each step of contraction only differs by one triangle in the process. The homeomorphic mapping from inside of $C$ that contains finite number of triangles. They have a sequence of
homeomorphic mappings from one to its neighbor in contraction (adding a triangle) . We will end up with a final triangle that is  homeomorphic to a disk.
So we complete the proof.

The key part of this proof is to construct an algorithm: We use the graph-distance for the construction of contraction sequence.

\begin{theorem}(Simple Connectedness of the Euclidean Plane)
 The Euclidean Plane is simply connected.
\end{theorem}

\begin{proof}

(According to the definition of the simple connectedness, we want to prove that each simple closed curve can be contracted to be a point (0-cell). The following proof is to 
describe an algorithm to perform such a contraction. In fact, when our contraction algorithm reaches the boundary of a triangle, we can see that we just need to shrink this 
triangle to be a point (0-cell). Since we already proved the Jordan theorem: A simple closed curve on two dimensional plane separates the plane into two connected components
in discrete case, we only need to guarantee that: (1) A 2D plane can be triangulated. This is true. (2) Any simple closed cycle can be contracted to a point. The contraction 
process is formed by a sequence of consecutive discrete curves. Any two adjacent curves only differs by a triangle or a set of triangles as long as these two curves are gradually
varied. The triangle (or triangles) are part of the initial triangulation. The meaning of ``differs'' here is $modulo2$ of the two curves is the boundary of a triangle.)

 In this proof, we first construct a triangulation that will made the original simple polygon $C$ as the boundary of the triangulation. There are many ways to make the triangulation in
 computational geometry with adding new points inside of a polygon. We give a simple way here that is not difficult to implement.

 To be exact, we will actually use a huge square $S$ that contains $C$ . We request from any point of $C$ to the
 edge of the square is longer than the diameter of  $C$ ( the diameter of $C$ is the largest distance of two points in $C$). We only triangulate this big square in a very fine way such that
 each triangle's internal area at most contains one vertex point of $C$ inside, and each edge except two ending points at most contain one vertex point of $C$. To do this we will calculate
 the Euclidean distance of each pair of vertex points in $C$, find the minimum distance among these pairs $d_0$. We then made the edge of a regular triangle (for fine triangulation) as at most
 $1/3 \times d_0$ in its length.

 Since $C$ is a simple path, so
 there is no two vertices in $C$ they meet at one point in the triangulation. We also denote the square with the the triangulation as $S$.

 We now construct a special modification of the triangulation to make $C$ on the vertices of triangulations. Then we will made every angle of $C$ to be ``wide'' angle by inserting
 a triangle as necessary in order to satisfy the condition of Theorem 1.

 Step 1: If there is a vertex point of $C$ on the internal part of an edge (not at the ending points of an edge), we will make the intercept point to an actual point
 in triangulation of $S$. We link a line from this point to the third point (not on the intercepted edge) in two triangles  in $S$ that share the the intercepted edge.

 Step 2: If there is a vertex point of $C$ on the internal part of triangle , link this vertex to all three vertex points of the triangle.

 Step 3: If there is an edge of $C$ whose internal (not at the ending points) contains a vertex point in $S$, make this point as a vertex point in $C$. (Add a new point to $C$.)

 Step 4: If there is an edge of $C$ whose internal (not at the ending points) intersects with the internal part of an edge of a triangle in $S$, link the intersecting point
         to the third point (not on the intersecting edge) in two triangles  in $S$ that share the the intersecting edge. Make this intersecting point as a vertex point in $C$. (Add a new point to $C$.)

 Step 5: Repeat above steps to remove all cases mentioned in Step 1 to Step 4.

 The correctness of those steps are not hard to prove. In order to satisfy the conditions of Theorem 1, we refine the all triangles by adding a point in the central point to split a triangle into
 three triangles. If we do one more time, we will make each original angle in $C$ to be wide angle.  Now we can apply the easy version of the Jordan curve theorem (Theorem 1).

 According to the Jordan curve theorem, there will be two connected components in  $S-C$.

 We will first mark all points inside of $C$ by identify the
 connected component that contains the smaller (finite) numbers of triangles in the component. This also can be done by identifying a triangle vertex that
 does not connect to the boundary point in $S$.

 For a certain point $p$, we like to design the following algorithm to contract $C$ ($C$ might be modified to have more vertices). Calculate the distance from
 $p$ to all points on $C$ using only marked vertex points (The points only inside of $C$ now, we can declare.) All marked points and vertex points on $C$ (the new $C$) will
 make a vertex set $V_c$, all edges that has two marked ending points or on $C$ are collected as the edge set $E_c$. ($C$ will be a discrete curve since the ``wide'' angle property and
 three time of the length of edges between two vertex points in $C$.)  So we have $G_c =(V_c, E_c)$.

 The contraction will be made on graph $G_c =(V_c, E_c)$. The key idea is to use the graph-distance to find the furthermost 2-cell $\Delta$ (here is the triangle) on the planar graph---this
 is because $G_c$ is already embedded in $S$ that is a planar graph. As a planar graph, $G_c$ has the boundary of $C$. Then we will deduct this 2-cell $\Delta$. This deduction is special so that
 the new boundary after deduction will be a simple path $C_{1}$, i.e., $XorSum(C, C_1) = \Delta$. Therefore, $C$ and $C_1$ are gradually varied. We repeat this process, we will get to the single
 cell that contains $p$. The following construction will make  that $C_1$ is constructible using algorithmic technology.

 In the beginning, $C$ is a simple path.  In fact, $C$ is a discrete curve, more stronger than a simple path.  Find the graph-distances (the shortest length of edges in paths between two vertices)
 from $p$ to all points in $C$ in $G_c$. (We only care about this closed path. We no longer care too much about $S$.)  Now, there must be a point $x$ having the largest graph-distance to $p$. $x$ is contained in a triangle or several triangles in $G_c$.

 There are few cases: (1) $x$ is the only farthest vertex in $C$. (2) $x$ is one of the farthest vertices in $C$. Case (2) has two subcases: (i) The 2-cell containing $x$ does not contain any other $y$ that
 has the same distance to $p$. (ii) The 2-cell containing $x$  contains another $y$ that has the same distance to $p$.

 Note: Delete $x$ will not change the distance from $p$ to other vertex in $C$. this is because, there must be a shortest path to other points not passing $x$. The shortest path to other points passing $x$ will make
 that $x$ is not the farthest point in $C$.

 Case 1: $x$ is the only furthest vertex in $C$. (i) If $x$ is the vertex that is only contained in one triangle in $G_c$ ($C$ is the boundary), we can delete two edges linking to $x$ (of this triangle). Use the
 third edge to replace the two edges in $C$. So we get $C_1$ that is gradually varied to $C$. (ii)  If $x$ is the vertex that is contained in several triangles in $G_c$, we remove an edge containing $x$ in $C$, replaced by two other edges of the sample
 triangle that contains the removed edge. So the new path $C_1$ is gradually varied to $C$. And a triangle was removed from $G_c$, we continue this process until (i) occurs. Then we use the action in (i) to remove the point
 $x$. (One can also remove $x$ and all edges linking to $x$ use the half umbrella edges to replace the two edges on $C$ that contain $x$.)

 Case 2: Assume $x$ is the first such a point from clockwise of $p$,  (i) If $x$ is contained by a single triangle in $G_c$, we can delete two edges linking to $x$ (of this triangle). Use the
 third edge to replace the two edges in $C$. So we get $C_1$. ~\footnote  {A pathological case was found when we deal with a thin 2-manifold where each 1-cell on the boundary in a triangle having the following property: the third vertex of the triangle is also on the boundary.
 In such a case, we can not simply remove this triangle. However, we can always find a triangle having two 1-cells on the boundary. We need to remove that triangle first. See Appendix C in this paper or Appendix A in \cite{Che17}. This case only occurs where each triangle intersecting with the boundary curve $C_i$ for some $i$, the triangle containing the intersecting edge (1-cell) has the third point is also on $C_i$. We also noticed that
  using 1-cell distance would be better than using graph-distance that is 0-cell distance in the paper {\it L. Chen and S. Krantz,  A Discrete Proof of The General Jordan-Schoenflies Theorem,2015,}
\newline {\tt http://arxiv.org/abs/1504.05263}}
 (ii) If $x$ is the vertex that is contained in several triangles in $G_c$, we will first remove an edge if there is a neighbor in $C$ that has the same distance
 to $p$ as $x$. Otherwise, just remove an edge containing $x$ in $C$.  Replaced this removed edge by two other edges of the sample
 triangle that contains the removed edge. So the new path $C_1$ is gradually varied to $C$. And a triangle was removed from $G_c$, we continue this process until (i) occurs. Then we use the action in (i) to remove the point
 $x$. Remember any of removed edge will not be affect to the shortest path from $p$ to other points on $C$.

 We repeat the above process we will delete 2-cell one-by-one and get a sequence of gradually varied simple paths $C_1,C_2, \cdots, C_k$. Each of those paths will not cross-over any other.
 Since we only have finite number of 2-cells in $G_c$, the above process will be halt to end at the single 2-cell with the boundary $C_k$ containing $p$.

 So we proved that Euclidean plane is simply connected under our definition of discrete deformation. \qed

\end{proof}

A triangle is homeomorphic to a disk. So attach a triangle on an edge with sharing two vertices points will be homeomorphic to the first triangle. The homeomorphic mapping is
  done by dragging the middle point of the shared edge to the third point of the second triangle. Use the same procedure, we can attach another triangle to existing two. The  homeomorphic mapping is
   the same as the one by only dragging a shared edge. Therefore,  $C_k,C_{k-1}, \cdots, C_0=C$ defines the sequence of such homeomorphic mappings from the area bounded by $C_{i-1}$ to the area bounded $C_{i}$ in $G_c$. (We can define $G_{i}$ as the one where the edge and points were deleted. ) Thus, we have a sequence of invertible continuous mappings from $C$ to a disk. We have the following theorem:

\begin{theorem} ({\bf The Jordan–Schoenflies Theorem})   A simple closed curve on the Euclidean plane separate the plane into two connected components. The component with bounded closure
is homeomorphic to the disk.
\end{theorem}

Please also note that the proof of simple connectedness of the Euclidean plane may already done by others in some discrete way. To put a proof in this paper as an appendix is to make
the proof of original Jordan curve theorem self-contained in this paper.

To consider a proof of the Jordan–Schoenflies theorem was inspirited by the discussion with Professor F. Luo and X. Huang.  Professor Steven G. Krantz also mentioned the author this theorem through email communications. Moise
had a similar proof of the Jordan–Schoenflies theorem in \cite{Moi}. \\

\section {Appendix  C:  A Special Case in the Proof of The Jordan–Schoenflies Theorem}

In \cite{Che17}, when we deal with high dimensional contraction, we found a special case that should be discussed when we do the contraction in the proof of Theorem 3 in Appendix B. Now we just refer to the context in \cite{Che17}:

{\bf
Before we prove Theorem A for 3D cases or higher dimensional cases, we first prove it for 2D cases. Because in 2D, Case 2 does not exist (Otherwise we admits two 1-cells on $B$ already).

{\bf Theorem B:} There is an $2$-cell of $D$ having two $1$-cell on Boundary $B$ of $D$ if for each 2-cell $e$ having a 1-cell (2-face) $f$ on $B$, the following condition holds:
: The 0-cell $P\in e$ not in $f$ is on $B$ .

{\it Proof:}

We define the distance between two cells $E$ and $E'$ is the shortest distance (the length of a shortest edge-path or 1-cell path) between two points in each $E$ and $E'$. $E$ and $E'$ can be in different dimensions.  Denote
\newline $d_{\Omega}(E,E')=\{ n | n $  is the length of a shortest path from a 0-cell \newline\hskip 2in in $E$ to another 0-cell in $E' \in \Omega \}$.

We now recall the concept of $i$-cell-distance: the length of the shortest $i$-cell path between two cells (usually 0-cells, include these two cells, these two cells are not the same) where each adjacent pair of two $i$-cells share an $(i-1)$-cell. In general we use $d_{\Omega}^{(i)}(E,E')$ for $i$-cell length that always shares an $(i-1)$-cell between $i$-cells in the path.  Use $d_{D}^{(i)}(X_0,X_1)$ to denote this distance for two points or two cells.  For instance, 1-cell distance is the graph distance. We know that each pair of 3-cells in $D$ is 2-connected.
We can identify two points on $B$ that has the largest 3-cell-distance in $D$.
We will prove that we can find $e$ nearby one of the two points that has two 2-faces on $B$.

We start with the example in two-dimensional cases. See Fig.~\ref{fig:TwoDCase}

This manifold has the following property: (a) there is no inner point, and
(b) each 1-cell on boundary has associated 2-cell that contains this 1-cell. This 2-cell has the third point on the boundary too.

\begin{figure}[h]
\begin{center}
   \epsfxsize=4.5in
  \epsfbox{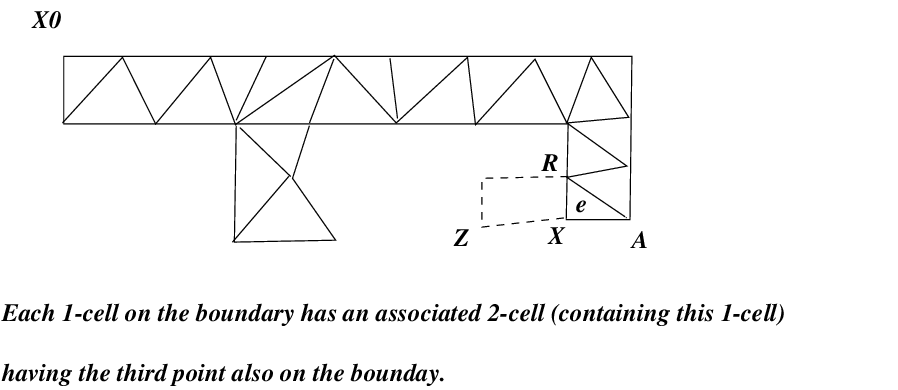}   
\caption{ A 2D example where exists a 2-cell that must have two 1-cells on the boundary $B$.} \label{fig:TwoDCase}
\end{center}
\end{figure}

{\bf We now to prove for the 2D case: the first proof}
The assumption is: For each 2-cell $e \in D$, if $e$ has a 1-cell $f$ in $B$ ($B$ is a 1-cycle) and the third point $P$ of $e$ not in $f$ is also on $B$, then there must be a 2-cell $e'$ that has  has two 1-cells on $B$.
So, we can do the contractional removal of the corresponding 2-cell $e'$ to maintain that the new boundary is still a 1-cycle after removal of $e'$ .

We will use 2-cell distance. We will see that the far-most pair of points in $D$ in 2-cell distance will help to determine such an $m$-cell, $m=2$, that has two $(m-1)$-cells on $B$. The good thing in 2D is that a 2-cell, having a 1-cell on $B$ in $D$, contains a 1-cell not on $B$ that will split $D$ into two disconnected parts. Assume that these two points are $X_0$ and $X$ (meaning that this pair has the longest 2-cell distance). Let the 2-cell containing $X$ is $e=ARX$. If $e$ does not have two 1-cells on $B$, we can assume the 1-cell $AX$ is on $B$ ($R$ is at another side of $B$ from $A$). See Fig. ~\ref{fig:TwoDCase}. Since $B$ is a 1-cycle, then there must exist another 1-cell on $B$ that contains $X$. (Any 0-cell is contained by exact two 1-cells.) This 1-cell can be named $XZ$, $Z$ is on $B$, and from $Z$ to $R$ there is a 1-cell path (since $B$ is connected). $Z$ is the point that has strictly longer 2-cell-distance comparing to $X$ from $X_0$ since any 2-cell path to $Z$ from $X_0$ must include cell $e$. So, we get the contradiction that $X$ is the far-most point from $X_0$. Therefore, there must be a 2-cell that has two 1-cells on $B$.

Thus, we proved the 2D case. The idea of proving for the 3D case can be similar as the case for 2D but having more complexity.  The above simple proof for 2D can be inserted to \cite{Chen16, Chen14} as a supplement fact for the 2D case.  

Note that in this proof, we can see that we only need to get a pair that has the local maximum distance (in $Star_D(X)$ from $X_0$) in 2-cell path of $D$. In addition, $X_0$ can be any point on $B$.

(We also have another proof of this theorem. See the Appendix A of the related paper \cite{Che17} where the revision number was v5. We omitted here.) \qed}


\end{document}